\documentclass[11pt,a4paper,reqno]{amsart}
\usepackage{graphicx} 
\usepackage{amsmath,amssymb,amsthm}
\usepackage{xcolor}
\usepackage{todonotes}\usepackage[normalem]{ulem}
\usepackage{cite}
\usepackage{graphics}
\usepackage{appendix}
\usepackage[foot]{amsaddr}
\usepackage{algorithm}
\usepackage{algpseudocode}

\newtheorem{theorem}{Theorem}
\newtheorem{remark}[theorem]{Remark}
\newtheorem{lemma}[theorem]{Lemma}
\newtheorem{proposition}[theorem]{Proposition}

\newtheorem{corollary}[theorem]{Corollary}
\usepackage{cancel}

\newcommand{\changed}[1]{\textcolor{black}{#1}}
\newcommand{\changedsecond}[1]{\textcolor{black}{#1}}
\usepackage[european]{circuitikz}
\usetikzlibrary{decorations.pathreplacing}

\usepackage{geometry}
\numberwithin{equation}{section}
\allowdisplaybreaks

\newcommand{\R}{\mathbb{R}}
\newcommand{\N}{\mathbb{N}}

\title{Goal-oriented time adaptivity for port-Hamiltonian systems$^*$}

\author[A.\ Bartel]{Andreas Bartel}
\address{{\bfseries A.~Bartel:} Institute of Mathematical Modelling, Analysis and Computational Mathematics (IMACM), School of Mathematics, Bergische Universität Wuppertal, Wuppertal, Germany}

\author[M.\ Schaller]{Manuel Schaller}
\address{{\bfseries M.~Schaller:} Institute of Mathematics, Technische Universit\"at Ilmenau, Ilmenau, Germany}
\thanks{$^*$This work was supported within the CRC 1701 "Port-Hamiltonian Systems".}

\begin{document}

\maketitle

\begin{abstract}
     Port-Hamiltonian systems provide an energy-based modeling paradigm for dynamical input-state-output systems. At their core, they fulfill an energy balance relating stored, dissipated and supplied energy. To accurately resolve this energy balance in time discretizations, we propose an adaptive grid refinement technique based on a posteriori error estimation. 
The evaluation of the error estimator includes the computation of adjoint sensitivities. 
To interpret this adjoint equation as a backwards-in-time equation, we show piecewise weak differentiability of the dual variable. Then, leveraging dissipativity of the port-Hamiltonian dynamics, we present a parallelizable approximation of the underlying adjoint system in the spirit of a block-Jacobi method to efficiently compute error indicators. We illustrate the performance of the proposed scheme by means of numerical experiments showing that it yields a smaller violation of the energy balance when compared to uniform refinements and traditional step size controlled time stepping.
\end{abstract}



\section{Introduction}

\noindent In recent years, port-Hamiltonian systems (pHs) have received significant attention in various fields of mathematics and applications. At their core, pHs extend dissipative Hamiltonian equations to systems with inputs and outputs, so-called ports. These ports allow for control, observation and couplings of pHs by means of the exchanged power. 
Hence, pHs provide a modular framework particularly suitable for highly complex systems including different physical domains such as mechanics, electromagnetics, fluid dynamics and more. Moreover, due to their Hamiltonian nature, pHs satisfy an energy balance equation and are dissipative with respect to the power supply rate given by the product of the input and the conjugated output. This dissipativity property can be utilized in various disciplines of mathematics and control engineering such as functional analysis~\cite{JacoZwar12}, passivity-based controller \cite{OrteScha02} and observer design~\cite{VenkScha10} and (numerical) linear algebra \cite{MehlMehr18,MehrUnge23}. Furthermore, the energy balance allows to formulate energy-optimal control problems~\cite{SchaPhil21,FaulMasc22}, leading to highly efficient optimization-based controllers for pHs, e.g., for adaptive building control~\cite{SchaZell24}.

As pHs satisfy an energy balance in continuous time, it is desirable to resolve this energy balance also on a discrete level. For (closed) Hamiltonian systems, there is a rich literature on symplectic integration schemes achieving preservation of energy, cf.~\cite{HairHoch06,LeimReic04}. For pHs, which in particular include dissipative structures and external ports, there are recent advances using collocation methods \cite{MehrMora19,KotyLefe19}, see also the survey \cite{MehrUnge23} and the references therein. As the energy balance can usually not be fulfilled exactly on the discrete level, a common approach in discretizations is to require that it is not exact but consistent with the discretization scheme \cite{KotyLefe19} or to introduce discrete gradients \cite{GoreYalc08,KinoThom23,Gonz96}, where a surrogate energy balance holds for the time-discretized problem.

In this work, we present an approach to adaptively refine the time grid such that the energy balance should be approximately fulfilled without requiring a particular integrator. Our main motivation lies in the simulation of pHs by multirate schemes which are particularly designed to respect different time constants of subsystems which occur, e.g., in field/circuit coupling~\cite{SGB10} by means of individual step sizes~\cite{GeWe1984,BaGue2022}. Whereas first symplectic multirate schemes for systems without dissipation and ports are available \cite{SGS2023}, currently, it seems impossible for multirate schemes \cite{GueSan2016} to mimick an energy balance for pHs as in collocation or discrete gradient schemes.

Hence, instead of designing particular methods to respect the energy balance (by restricting the coefficients of a Runge-Kutta scheme or parameters of the collocation schemes), which might be too restrictive or even impossible for future applications in a pHs multirate context, we propose an approach using goal-oriented a posteriori adaptivity, as initiated in the seminal papers~\cite{Este95,Becker2000}. Since then, goal-oriented grid adaptivity has achieved wide success for ordinary and partial differential equations, e.g.\ in flow simulations \cite{BeckRann01,Hartmann2002}, optimal control and parameter estimation \cite{MeidVex07,Vexler2004,Kroener2011} and multirate methods~\cite{BaSc_Estep2012}. For a broad overview over the topic, we refer to the monograph \cite{Bangerth2003}. This method allows to perform grid refinement to obtain a small error of the solution in a user-specified functional, called the Quantity of Interest (QoI). The main contribution of this work is the application of goal-oriented techniques to linear finite-dimensional pHs by formulating a QoI targeting the violation of the energy balance.

The remainder of this manuscript is structured as follows. In Section~\ref{sec:ph_goal}, we provide the fundamentals of linear 
pHs and the main idea of goal-oriented adaptivity. Subsequently, in Section~\ref{sec:vardisc}, we deduce the variational formulation of the ODE and the corresponding discretization scheme by means of a Galerkin approach. Furthermore, we provide a result showing higher regularity of the adjoint state for right-hand sides of a particular structure. Then, we present a prototype implementation of goal-oriented methods for pHs in Section~\ref{sec:comp}. 
As this includes the computation of sensitivities by means of a linear system modeling a backwards-in-time equation, we propose in Section~\ref{sec:para} a method to approximate these sensitivities leveraging dissipativity. Last, we illustrate the performance of the proposed scheme by means of a numerical example in Section~\ref{sec:num}.

\section{Goal-oriented error estimation for port-Hamiltonian systems}\label{sec:ph_goal}
\noindent In this work, we consider linear port-Hamiltonian systems on finite time horizons $T>0$ with state $x:[0,T]\to \R^n$ and input $u:[0,T] \to \mathbb{R}^m$, $n,m\in \N$ of the form
\begin{subequations}\label{eq:ph}
	\begin{align}\label{eq:state}
	\dot{x}(t)&=(J-R)Qx(t) + Bu(t), \qquad x(0)=x^0,\\
	y(t) &= B^\top Qx(t). \label{eq:port}
	\end{align}
\end{subequations}
Here, $J=-J^\top\in \mathbb{R}^{n\times n}$ is the structure matrix, $0\leq R=R^\top\in \mathbb{R}^{n\times n}$ models the dissipation, $0\leq Q=Q^\top\in \mathbb{R}^{n\times n}$ corresponds to a quadratic
Hamiltonian $H(x)= \frac{1}{2} x^\top Q x$ and $B \in \mathbb{R}^{n\times m}$ is an input matrix. For given input $u\in L^2([0,T],\mathbb{R}^m)$ and initial value $x^0\in \R^n$, we consider the weak solution $x \in H^1([0,T],\mathbb{R}^n)$ such that \eqref{eq:ph} is satisfied for a.e.~$t\in [0,T]$. 

By the structure of the involved matrices and direct calculations, the solution of \eqref{eq:ph} obeys the power balance
\begin{align}\label{eq:power}
\frac{\mathrm{d}}{\mathrm{d}t} H(x(t)) = u(t)^\top y(t) - \|R^{1/2}Qx(t)\|^2
\end{align}
for 
a.e.~$t\in [0,T]$ and the corresponding energy balance
\begin{align}\label{eq:energy}
H(x(t)) - H(x(s)) = \int_s^t  u(\tau)^\top y(\tau) - \|R^{1/2}Qx(\tau)\|^2\,\mathrm{d}\tau
\end{align}
for all $s,t\in [0,T]$, where $R^{1/2}$ denotes the matrix square root of $R=R^\top \geq 0$.

In case of a closed system without dissipation, i.e., $R=0$ and $u=0$, where \eqref{eq:state} reduces to a conservative Hamiltonian system, then this balance can be exactly fulfilled on the discrete level when using a symplectic integration scheme, since $H$ is a quadratic invariant. More precisely, if $\tilde{x}: \mathcal{T} \to \R^n$ is the numerical approximation of the true state $x:[0,T]\to \R^n$ on a time grid $\mathcal{T}$ with elements $0=t_0<\ldots<t_M=T$, $M\in \N$, we have, using, e.g., and implicit midpoint rule for \eqref{eq:ph} with $R=0$ and $u=0$,
\begin{align}
H(\tilde{x}(t)) = H(\tilde{x}(0)) \qquad \forall t\in \mathcal{T}. 
\end{align}
Using the mature literature on symplectic integration, e.g.~\cite{HairHoch06,LeimReic04}, such a property can also be obtained in nonlinear cases. For systems involving dissipation and inputs, i.e., port-Hamiltonian systems, we refer to the work \cite{KotyLefe19} which formulates criteria for collocation-based discretizations to obtain an energy balance in discrete time. 

The focus of this work is complementary to these approaches. Instead of designing methods satisfying an energy balance, we suggest an approach to design adaptive grids via a posteriori error estimation that allow to resolves the above balance as well as possible. 

\subsection{Refinements targeting the energy balance}
\noindent In this part, we propose the main idea of using goal-oriented error estimation~\cite{Becker1996,Bangerth2003} for pHs. In the following, we denote by $x:[0,T]\to \R^n$ the solution of \eqref{eq:ph} and by $\tilde x :[0,T]\to \R^n$ a numerical approximation. Roughly speaking, goal-oriented adaptivity allows to adaptively refine the time grid such that the error of the numerical approximation measured in a real-valued Quantity of Interest (QoI), denoted by $I$, becomes smaller than a prescribed tolerance $\mathrm{tol}>0$, i.e.,
\begin{align}\label{eq:qoismall}
|I(\tilde x) - I(x)| \leq \mathrm{tol} .
\end{align}
For port-Hamiltonian systems, a first natural quantity of interest can be deduced from the energy balance \eqref{eq:energy}. 
Let $\tilde{y}=B^\top Q \tilde{x}$ be the output of the approximation $\tilde x$. Then, we define the QoI corresponding to the violation on the full time horizon by
\begin{align*}
I_\mathrm{glob}(\tilde{x}) := \left|\int_0^T \! -u(t)^\top \tilde{y}(t)  + \tilde x(t)^\top \changed{QRQ} \tilde x(t)\,\mathrm{d}t 
\;+ \;   H(\tilde{x}(T)) - H(\tilde{x}(0))\right|^2 .
\end{align*}
In the view of the energy balance \eqref{eq:energy}, the solution $x\in H^1([0,T],\R^n)$ of \eqref{eq:ph} satisfies $I_\mathrm{glob}(x)=0$. Hence, \eqref{eq:qoismall} reduces to the requirement that
\begin{align*}
|I_\mathrm{glob}(\tilde x)|\leq \mathrm{tol}.
\end{align*}
An alternative quantity is given by piecewise consideration of the energy balance.
For a suitably chosen a priori grid $0=t_0<\ldots,t_M=T$, we set
\begin{align}\label{eq:Iloc}
I^i_{\mathrm{loc}}(\tilde x) := 
\left|\int_{t_i}^{t_{i+1}} \hspace{-.5cm} -u(t)^\top \tilde y(t) + \tilde x(t)^\top \changed{QRQ}\tilde x(t)\,\mathrm{d}t 
+ H(\tilde x(t_{i+1})) \hspace{-.05cm} - \hspace{-.05cm}H(\tilde x(t_i))
\right|^2
\end{align}
for $i=0,1,\ldots,M-1$. Then, the corresponding total error reads
\begin{align} \label{eq:Iloctotal}
I_{\mathrm{loc}}(\tilde x) := \sum_{i=0}^{M-1} 
I^i_{\mathrm{loc}}(\tilde x). 
\end{align}
By the triangle inequality, we directly observe that if we control the total local error, we also control the global error whereas the converse is not true, as violations of the local energy balances might cancel out in $I_\mathrm{glob}$. 
Furthermore, we note that $I_{\mathrm{loc}}(x)$ mimics better the power balance in continuous time~\eqref{eq:power}, or, correspondingly the energy balance \eqref{eq:energy} on all subintervals. Thus, in the remainder, we will focus on the QoI given by $I_\mathrm{loc}$.

\noindent To ensure that not only the energy balance is satisfied, but to also \changedsecond{guarantee a small error in norm}, 
we further consider a weighted quantity of interest, i.e., for $\rho \ge 0$, we set
\begin{align}\label{eq:Ilocrho}
I_{\mathrm{loc},\rho}(x) := I_\mathrm{loc}(x) +  \rho \,\|x\|^2_{L^2([0,T],\R^n)}.
\end{align}
\changed{In view of \eqref{eq:qoismall}, this QoI additionally yields error control in view of the norm.}


\subsection{Fundamentals of goal-oriented adapativity}
\noindent To motivate the variational formulations in the upcoming section, we briefly introduce the fundamentals of goal-oriented adaptivity and refer the reader to \cite{Bangerth2003} or \cite{BeckRann01} for further details. Let $X$ be a reflexive Banach space and consider $\mathcal{A}:X\to X^*$ linear and bounded as well as $f\in X^*$ be given. Later, this operator $\mathcal{A}$ will model the state dependent part of the ODE~\eqref{eq:state} including time derivative and initial condition. Moreover, we consider a discrete space $X_h\subset X$ with a Galerkin projection of $\mathcal{A}$ given by $\mathcal{A}_h:X_h\to X_h^*$. Let $x\in X$ be the solution of a variational problem
\begin{align}\label{eq:varprobx}
\langle \mathcal{A}x,\varphi\rangle_{X^*,X} = \langle f,\varphi\rangle_{X^*,X} \quad \forall \varphi\in X
\end{align}
and $x_h\in X_h$ be the solution of the projected problem
\begin{align}\label{eq:varprobxproj}
\langle \mathcal{A}_h x_h,\varphi_h\rangle_{X_h^*,X_h}  
= \langle f_h,\varphi_h \rangle_{X_h^*,X_h} \quad \forall \varphi_h\in X_h .
\end{align}
Consider a Fréchet differentiable QoI $I:X\to \R$. We define an adjoint state $\lambda \in X$ solving
\begin{align}\label{eq:adjoint}
\langle \mathcal{A}^*\lambda,\varphi \rangle_{X^*,X} = \langle I'(x_h),\varphi\rangle_{X^*,X} \quad \forall \varphi\in X,
\end{align}
where we identified $X^{**}\cong X$ due to reflexivity, and an approximation $\lambda_h\in X_h$ solving
\begin{align}\label{eq:adjoint_projected}
\langle \mathcal{A}_h^*\lambda_h,\varphi_h \rangle_{X_h^*,X_h} = \langle I'(x_h),\varphi_h\rangle_{X_h^*,X_h} \quad \forall \varphi_h\in X_h .
\end{align}
Then, 
the error in the QoI may be approximated via
\begin{align}\label{eq:errorestimator}
I(x) - I(x_h) \approx \langle \mathcal{A}(x-x_h),\lambda \rangle_{X^*,X} =  \langle \mathcal{A}(x-x_h),\lambda-\lambda_h\rangle_{X^*,X},
\end{align}
cf.~\cite[Section 2.2]{BeckRann01}. Note that in the last equality we used Galerkin orthogonality, i.e., $\langle \mathcal{A} (x-x_h),\lambda_h\rangle_{X^*,X} = 0$. 

Thus, evaluation of the error estimator and computation of local refinement indicators amounts to two steps (for details we refer to Section~\ref{sec:comp}):
\begin{enumerate}
	\item[(i)] Dual problem: Solve the projected adjoint system~\eqref{eq:adjoint_projected} to obtain $\lambda_h$ and approximate $\lambda$ solving \eqref{eq:adjoint}, e.g., using a higher order method. 
	\item[(ii)] Evaluation and localization: Approximate $x$ solving \eqref{eq:varprobx} from $x_h$ solving \eqref{eq:varprobxproj} and evaluate the error \eqref{eq:errorestimator}. In many applications the duality product in \eqref{eq:errorestimator} reduces to a sum over the grid cells (e.g., time intervals) directly leading to local error indicators that can be used for adaptive refinements.
\end{enumerate}

In the following section, we provide a suitable variational formulation of \eqref{eq:ph} to apply the above goal-oriented methodology. 

\section{Variational formulation and time discretization}\label{sec:vardisc}
\noindent  In this section, we first provide the formulation of the state equation~\eqref{eq:state} in a variational framework in Section~\ref{subsec:state}. Then, in Section~\ref{subsec:adjoint}, we deduce the adjoint equation necessary to evaluate the sensitivities w.r.t.~the QoI by means of \eqref{eq:adjoint}. For QoIs with a particular structure, such as $I_\mathrm{loc}$ as defined in \eqref{eq:Iloc}, we show that the adjoint state enjoys higher regularity, i.e., it is piecewise weakly differentiable with jumps given by the point evaluations of the right-hand side. Last, in Section~\ref{subseq:disc},  we provide a time discretization by means of a non-conforming Galerkin approach as presented in \cite{MeidVex07} for parabolic partial differential equations.

\subsection{State equation}\label{subsec:state}
\noindent We formulate the dynamics \eqref{eq:state} in a variational setting to enable goal-oriented methods.  To this end, we abbreviate $A:=(J-R)Q$ and define the bounded and linear operator $\mathcal{A} : H^1([0,T],\R^n)\to L^2([0,T],\R^n)\times \R^n$ defined by
\begin{align*}
\langle \mathcal{A}x,(\varphi,\varphi_0)\rangle_{L^2([0,T],\R^n)\times \R^n} := \int_0^T \langle \dot x(t) - Ax(t), \varphi(t)\rangle\,\mathrm{d}t + \langle x(0),\varphi_0\rangle
\end{align*}
for all $(\varphi,\varphi_0)\in L^2([0,T],\R^n)\times \R^n$, where we tacitly identified $L^2([0,T],\R^n)\times \R^n$ with its dual space.
For input functions $u\in L^2([0,T],\R^m)$, we may thus formulate the dynamics~\eqref{eq:state} via the operator equation
\begin{align}\label{eq:odeoperator}
\langle \mathcal{A}x,(\varphi,\varphi_0)\rangle_{L^2([0,T],\R^n)\times \R^n} = \int_0^T \langle  Bu(t),\varphi(t)\rangle \,\mathrm{d}t + \langle x^0,\varphi_0\rangle
\end{align}
for all $(\varphi,\varphi_0)\in L^2([0,T],\R^n)\times \R^n$.
The following result shows that one can also formulate this variational problem with weakly differentiable test functions.
\begin{proposition}\label{prop:testfunc}
	If for all $\varphi\in H^1([0,T],\R^n)$,
	\begin{align*}
	\langle \mathcal{A}x,(\varphi,\varphi(0))\rangle_{L^2([0,T],\R^n)\times \R^n} = \langle f,\varphi\rangle_{L^2([0,T],\R^n)}+ \langle x^0,\varphi(0)\rangle.
	\end{align*}
	Then $\dot x = Ax+f$ a.e.\ on $[0,T]$ and $x(0)=x^0$.
\end{proposition}
\begin{proof}
	This directly follows from the density of $C_0^\infty([0,T];\R^n)$ in $L^2([0,T];\R^n)$ and the inclusion $C_0^\infty([0,T];\R^n)\subset H^1([0,T],\R^n)$.
\end{proof}

\subsection{Adjoint equation as piecewise backwards-in-time equation.}\label{subsec:adjoint}
\noindent  The adjoint operator $\mathcal{A}^*:L^2([0,T],\R^n)\times \R^n \to H^{1}([0,T],\R^n)^*$ is defined as
\begin{align*}
\langle \mathcal{A}^*(\lambda,\lambda_0),x\rangle_{H^1([0,T],\R^n)^*,H^1([0,T],\R^n)} = \langle \mathcal{A}x,(\lambda,\lambda_0)\rangle_{L^2([0,T],\R^n) \times \R^n}
\end{align*}
for $(\lambda,\lambda_0)\in L^2([0,T],\R^n)\times \R^n$ and for all $x\in H^1([0,T],\R^n)$. 

To interpret this adjoint equation~\eqref{eq:adjoint} as a backwards-in-time equation, we present a result proving higher regularity of the adjoint for particular right-hand sides consisting of $L^2$ terms and point evaluations. Such right-hand sides in particular occur when considering the derivative of $I_\mathrm{loc}$ as defined in~\eqref{eq:Iloctotal}. The presented analysis is strongly motivated by \cite{Schi2013}, where higher regularity of the adjoint was shown in the context of optimal control of parabolic partial differential equations. 

To define $I_\mathrm{loc}$ (cf.~\eqref{eq:Iloctotal}), we consider a time grid:
given $M\in \N$
\begin{subequations}\label{eq:time-discretization}
	\begin{align}
	\mathcal{T} &= \{t_0,t_1,\ldots,t_M\} \quad \text{ with  } \quad 0=t_0<t_1 < \ldots<t_M=T
	\intertext{with corresponding grid cells}
	\quad \mathbb{I}_i &= (t_{i-1},t_{i}] 
	\quad \text{for } \quad 1\leq i\leq M. 
	\end{align}
\end{subequations}
We introduce the respective space of piecewise differentiable functions: 
\begin{align*}
H^1_\mathcal{T}([0,T],\R^n) \!:=\! 
\{  v\in L^2([0,T],\R^n)  \,  \bigl|  \,v\vert_{\mathbb{I}_i}\in H^1( \mathbb{I}_i,\R^n)\ 
\forall 1\leq    i\leq M, v(0)\in \R^n
\}.
\end{align*}
For $\lambda \in H^1_\mathcal{T}([0,T], \R^n)$, we define the right and left limits at the points of discontinuity $t_i$ as follows: 
\begin{align*}
\lambda^-(t_i) &:= \lim_{t\to 0^+} \lambda(t_i - t)
\quad \text{(for }\, 0<i\le M), \\ \lambda^+(t_i) &:= \lim_{t\to 0^+} \lambda(t_i + t)
\quad \text{(for }\, 0\le i<M)
\end{align*}
and we set $\lambda^-(0):=\lambda(0)$. Furthermore, the jump at time instance $t_i$, $0\leq i \leq M-1$ is defined by
\begin{align*}
[\lambda]_i = \lambda^+(t_i) - \lambda^-(t_i).
\end{align*}
The following result shows that for right-hand sides consisting of an $L^2$-term and point evaluations, the adjoint state is piecewise weakly differentiable, i.e., $\lambda\in H^1_\mathcal{T}([0,T],\R^n)$ and it provides a characterization of the jumps at the points of discontinuity.
\begin{lemma}\label{lem:higherreg}
	Assume that there is $l\in L^2([0,T],\R^n)$ and $l_0,\ldots,l_M \in \R^n$ such that $(\lambda,\lambda_0)\in L^2([0,T],\R^n)\times \R^n$ satisfies
	\begin{align}\label{eq:backwards}
	\langle \mathcal{A}^*(\lambda,\lambda_0), \varphi\rangle_{H^1([0,T],\R^n)^*,H^1([0,T],\R^n)} = \langle l,\varphi\rangle_{L^2([0,T],\R^n)} + \sum_{i=0}^M l_i^\top \varphi(t_i) 
	\end{align}
	for all $\varphi \in H^1([0,T],\R^n)$. Then 
	\begin{enumerate}
		\item $\lambda\in H^1_\mathcal{T}([0,T],\R^n)$ with $-\dot \lambda - A^\top \lambda = l$ on $\mathbb{I}_i$ for all $1\leq i\leq  M$,
		\item $\lambda^+(0) = \lambda_0 - l_0$, $\lambda^-(T) = l_M$, and
		\item $[\lambda]_i = -l_i$ for all $0\leq i\leq  M-1$. 
	\end{enumerate} 
\end{lemma}
\begin{proof}
	We proceed analogously to \cite[Proposition 2.5, Proposition 3.8]{Schi2013}. To this end, let $1\leq i \leq M$ and $\varphi \in C^\infty([0,T],\R^n)$ such that $\varphi(t) =0$ for all $t\notin (t_{i-1},t_{i})$. Then
	\begin{align*}
	\langle l,\varphi\rangle_{L^2([t_{i-1},t_{i}],\R^n)} &= \langle l,\varphi\rangle_{L^2([0,T],\R^n)} + \sum_{i=0}^M l_i^\top \varphi(t_i) \\
	&= \langle \mathcal{A}^*(\lambda,\lambda_0), \varphi\rangle_{H^1([0,T],\R^n)^*,H^1([0,T],\R^n)} \\
	&= \langle \mathcal{A}\varphi,(\lambda,\lambda_0)\rangle_{L^2([0,T],\R^n)\times \R^n} = \langle \dot \varphi - A\varphi,\lambda\rangle_{L^2([t_{i-1},t_{i}],\R^n)}.
	\end{align*}
	This implies that, for all $1\leq i\leq M$, $-(l+A^\top \lambda)\vert_{\mathbb{I}_i}\in L^2([t_{i-1},t_{i}],\R^n)$ is the weak derivative of $\lambda\vert_{\mathbb{I}_i}$, i.e., in particular $\lambda\in H_\mathcal{T}^1([0,T],\R^n)$. Hence, we may perform integration by parts on the subintervals:
	\begin{align*}
	\langle \mathcal{A}^*&(\lambda,\lambda_0), \varphi\rangle_{H^1([0,T],\R^n)^*,H^1([0,T],\R^n)} 
	\\
	& = \langle \mathcal{A}\varphi,(\lambda,\lambda_0)\rangle_{L^2([0,T],\R^n)\times \R^n}
	\\
	& = \int_{0}^T \langle \dot \varphi(t)-A\varphi(t),\lambda(t)\rangle\,\mathrm{d}t + \langle \varphi(0),\lambda_0\rangle
	\\
	& = \sum_{i=1}^M \int_{\mathbb{I}_i} \langle \dot \varphi(t),\lambda(t)\rangle \,\mathrm{d}t + \langle \varphi(0),\lambda_0\rangle - \int_0^T \langle A\varphi(t),\lambda(t)\rangle\,\mathrm{d}t
	\\
	& = -\sum_{i=1}^M \int_{\mathbb{I}_i}  \langle  \varphi(t),\dot \lambda(t)\rangle \,\mathrm{d}t 
	+ \sum_{i=1}^{M} \left( 
	\langle \varphi(t_{i}),\lambda^-(t_{i})\rangle
	- \langle \varphi(t_{i-1}),\lambda^+(t_{i-1})\rangle 
	\right)\\
	&\qquad \qquad +  \langle \varphi(0),\lambda_0\rangle- \int_0^T \langle A\varphi(t),\lambda(t)\rangle\,\mathrm{d}t.
	\end{align*}
	For the jump terms we compute
	\begin{align*}
	\sum_{i=1}^{M} \bigl( 
	\langle &\varphi(t_{i}),\lambda^-(t_{i})\rangle 
	- \langle \varphi(t_{i-1}),\lambda^+(t_{i-1})\rangle
	\bigr)
	+ \langle \varphi(0),\lambda_0\rangle\\
	&=\sum_{i=1}^{M-1} \langle \varphi(t_{i}),\lambda^-(t_{i})-\lambda^+(t_i)\rangle + \langle \varphi(T),\lambda^-(T)\rangle + \langle \varphi(0),\lambda_0-\lambda^+(0)\rangle .
	\end{align*}
	Then, the result follows from \eqref{eq:backwards} as
	\begin{align*}
	E : H^1([0,T],\R^n) &\to L^2([0,T],\R^n) \times \R^n\times \ldots\times \R^n\\
	\varphi&\mapsto (\varphi,\varphi(t_0),\ldots,\varphi(t_M))
	\end{align*}
	has dense range, cf.~\cite[Lemma 2.3]{Schi2013}.
\end{proof}
Note that the previous result shows that, whereas the adjoint is given by a tuple $(\lambda,\lambda_0)\in L^2([0,T],\R^n)\times \R^n$, we may recover $\lambda_0$ from the right-hand side and the solution $\lambda$ via $\lambda_0 = \lambda^+(0) + l_0$.

Motivated by the higher regularity result of Lemma~\ref{lem:higherreg}, we may define a backwards-in-time operator $\mathcal{A}^- : H^1_\mathcal{T}( [0,T],\R^n) \to L^2([0,T],\R^n)\times \R^n \times \ldots \times \R^n$ via
\begin{align*}
\langle \mathcal{A}^- &\lambda,(\varphi,\varphi_0,\ldots,\varphi_T)\rangle_{L^2([0,T],\R^n)\times \R^n \times \ldots \times \R^n} 
\\
&:=  \sum_{i=1}^M \int_{\mathbb{I}_i} \langle -\dot \lambda(t) - A^\top \lambda(t),\varphi(t)\rangle\,\mathrm{d}t - \sum_{i=1}^{M-1} \langle [\lambda_i],\varphi_i\rangle + \langle \lambda(T),\varphi(T)\rangle.
\end{align*}
Directly from the proof of Lemma~\ref{lem:higherreg}, we obtain the following counterpart of Proposition~\ref{prop:testfunc} for the adjoint state.
\begin{corollary}
	If for all $\varphi\in H^1([0,T],\R^n)$,
	\begin{align*}
	\langle \mathcal{A}^-\lambda,(\varphi,\varphi(0),\ldots,\varphi(T))\rangle_{L^2([0,T],\R^n)\times \R^n \times \ldots\times \R^n} = \langle f,\varphi\rangle_{L^2([0,T],\R^n)}+  \sum_{i=0}^M l_i^\top \varphi(t_i).
	\end{align*}
	Then $-\dot \lambda = A^\top\lambda+f$ on $\bigcup_{i=0}^{M-1}(t_i,t_{i+1})$, $\lambda(T) = l_M$ and $[\lambda_i] = -l_i$ for $0 \leq i \leq M-1$. Moreover, for all $\lambda\in H^1_\mathcal{T}([0,T],\R^n)$ and $x \in H^1([0,T],\R^n)$ we have
	\begin{align*}
	\langle A^*(\lambda,\lambda(0)),x&\rangle_{H^1_\mathcal{T}([0,T],\R^n)^*,H^1_\mathcal{T}([0,T],\R^n)} \\&= \langle \mathcal{A}^- \lambda, (x,x(0),\ldots,x(T))\rangle_{L^2(0,T;\R^n)\times \R^n  \times \ldots \times \R^n} .
	\end{align*}
\end{corollary}

\subsection{Discretization}\label{subseq:disc}
\noindent In the following, we focus on discretization by the implicit Euler method, which can be viewed as a Galerkin projection onto a basis of piecewise constant functions. As these functions are discontinuous, they are not contained in the solution space $H^1([0,T],\R^n)$ such that this resembles a non-conforming Galerkin scheme, cf.~\cite{MeidVex07}.

To formulate the dynamics \eqref{eq:state} (or equivalently \eqref{eq:odeoperator}) for functions with discontinuities, we provide first an alternative formulation of the operator $\mathcal{A}$. To this end, we apply integration by parts such that for all $\varphi\in H^1([0,T],\R^n)$, 
\begin{align}
\langle \mathcal{A}x,&(\varphi,\varphi(0))\rangle_{L^2([0,T],\R^n)\times \R^n}\nonumber\\&= \int_0^T \langle \dot x(t) - Ax(t), \varphi(t)\rangle \,\mathrm{d}t + \langle x(0),\varphi(0)\rangle \nonumber\\
&= \langle \dot x, \varphi\rangle_{L^2([0,T],\R^n)}+ \langle x(0),\varphi(0)\rangle - \langle Ax,\varphi \rangle_{L^2([0,T],\R^n)}\nonumber\\
&= -\langle x,\dot \varphi \rangle_{L^2([0,T],\R^n)} + \langle x(T),\varphi(T)\rangle - \langle Ax,\varphi \rangle_{L^2([0,T],\R^n)}.\label{eq:startingpoint}
\end{align}
This formula now may be used to define the discretized version of $\mathcal{A}$, as it does not hinge on differentiability of the state $x$. We proceed similarly to the previous Section~\ref{subsec:adjoint}
and introduce a time grid $\mathcal{T}$ with cells $\mathbb{I}_i$ as in \eqref{eq:time-discretization} and step sizes $h_i = t_i - t_{i-1}$ ($1\leq i\leq M$).
We define the space of ansatz and test functions for the Galerkin ansatz by
\begin{align}\label{eq:Xh}
X_h \!:=\! \big\{ v_h \in L^2([0,T],\R^n) \,|\, v_h \vert_{\mathbb{I}_i}  =v_{h,i} \in \R^n, 1\leq i \leq M, \ v_{h}(0) \in \R^n\big\}.
\end{align}
Note that $X_h \subset H^1_\mathcal{T}([0,T],\R^n)$ as introduced in Section~\ref{subsec:adjoint}, as piecewise constant functions are clearly piecewise weakly differentiable. Moreover, $X_h$ is finite dimensional and we may identify
\begin{align*}
X_h \simeq \R^{n} \times \ldots \times \R^n = (\R^n)^{M+1} .
\end{align*}
Thus, for $x_h\in X_h$ and all $0 \leq i \leq M-1$, we have
\begin{align*}
x_h^-(t_i) = \lim_{t \to0^+}x_h(t_i-t) = \changed{x_h}(t_i),\qquad x_h^+(t_i) = \lim_{t\to 0^+} x_h(t_i+t) = \changed{x_h}(t_{i+1})
\end{align*}
and we define the jump at time instance $t_i$, $i=0,\ldots,M-1$ by 
\begin{align*}
[x_h]_i := x_h^+(t_{i}) -  x_h^-(t_{i}) = \changed{x_h}(t_{i+1}) - \changed{x_h}(t_i).
\end{align*}

\noindent \textbf{Discretized state equation.}  In view of \eqref{eq:startingpoint}, we compute for all
$x_h\in X_h$ and $\varphi \in H^1([0,T],\R^n)$
\begin{align*}
&-\langle x_h,\dot \varphi \rangle_{L^2([0,T],\R^n)} + \langle x_h(T),\varphi(T)\rangle\\
& \qquad = - \sum_{i=1}^M \int_{\mathbb{I}_i} \langle x_h(t),\dot \varphi(t)\rangle \,\mathrm{d}t + \langle x_h(T),\varphi(T)\rangle\\
& \qquad = \sum_{i=1}^M \int_{\mathbb{I}_i} \langle \dot x_h(t),\varphi(t)\rangle \,\mathrm{d}t + \sum_{i=1}^M \left(-\langle x_h^-(t_i),\varphi(t_i)\rangle + \langle x_h^+(t_{i-1}),\varphi(t_{i-1})\rangle \right) 
\\& \hspace{2cm}+ \langle x_h(T),\varphi(T)\rangle\\
& \qquad = \sum_{i=1}^{M-1} \langle x_h^+(t_{i}),\varphi(t_{i})\rangle -  \langle x_h^-(t_{i}),\varphi(t_{i})\rangle + \langle x_h^+(0),\varphi^+(0)\rangle\\
& \qquad = \sum_{i=1}^{M-1} \langle x_h^+(t_{i})-  x_h^-(t_{i}),\varphi(t_{i})\rangle + \langle x_h^+(0),\varphi(0)\rangle\\
& \qquad = \sum_{i=0}^{M-1} \langle x_h^+(t_{i})-  x_h^-(t_{i}),\varphi(t_{i})\rangle + \langle x_h^-(0),\varphi(0)\rangle\\   
& \qquad = \sum_{i=0}^{M-1} \langle [x_h]_i,\varphi(t_{i})\rangle + \langle x_h^-(0),\varphi(0)\rangle. 
\end{align*}
Thus, we define the Galerkin projection of $\mathcal{A}$ denoted by $\mathcal{A}_h : X_h \to X_h^*$ via
\begin{align*}
\langle \mathcal{A}_hx_{h},\varphi_h\rangle_{X_h^*,X_h} := \sum_{i=1}^M h_i \langle -Ax_h(t_i), \varphi_h(t_i) \rangle &+  \sum_{i=0}^{M-1} \langle [x_h]_i,\varphi_h(t_{i+1})\rangle \\&+ \langle x_h(0), \varphi_h(0)\rangle.
\end{align*}
for all $x_h, \varphi_h\in X_h$.
Hence, the Galerkin projection of the ODE \eqref{eq:state} (or equivalently \eqref{eq:odeoperator}) onto $X_h$ is given by
\begin{align*}
\sum_{i=1}^M h_i\langle- Ax_{h}(t_i), \varphi_h(t_i) \rangle\,\mathrm{d}t &+ \sum_{i=0}^{M-1} \langle [x_h]_i,\varphi_h(t_{i+1})\rangle + \langle x_h(0), \varphi_h(0)\rangle \\
&=  \sum_{i=1}^M \int_{\mathbb{I}_i} \langle Bu, \varphi_h(t_i)\rangle\,\mathrm{d}t + \langle x^0,\varphi_h(t_0)\rangle \quad \forall \varphi \in X_h.
\end{align*}
By choosing test functions satisfying $\varphi_h(t_j) = 0$ for all $j\neq i$, we can deduce decoupled equations for all time instances $t_i$, $1\leq i\leq M$. In fact, when combined with a rectangular rule for the integral including the input, this gives the implicit Euler method. More precisely, solving the above discretized ODE amounts to solving the linear equation system
\begin{align}\label{eq:state_disc}
\left(\begin{smallmatrix}
I & 0 & 0 & 0 & \cdots & 0 \\[0.5ex]
-I & I-h_1A & 0 & 0 & \cdots & 0 \\[0.5ex]
0 & -I & I-h_2A & 0 & \cdots & 0 \\[0.25ex]
\vdots & \vdots & \ddots & \ddots & \ddots & \vdots \\[0.25ex]
\vdots & \vdots & \vdots & -I & I-h_{M-1} A & 0 \\[0.25ex]
0 & 0 & 0 & 0 & -I & I-h_{M} A  \\
\end{smallmatrix}\right)
\left(\begin{smallmatrix}
x(t_0)\\[0.25ex] 
x(t_1)\\[0.25ex] 
x(t_2) \\[0.25ex]
\vdots\\[0.5ex] 
x(t_{M-1})\\[0.25ex]
x(t_M)
\end{smallmatrix}\right)
= 
\left(\begin{smallmatrix}
x^0\\[0.25ex] 
h_1 Bu(t_1)\\[0.25ex] 
h_2 Bu(t_2) \\[0.25ex]
\vdots\\[0.5ex] 
h_{M-1} Bu(t_{M-1})\\[0.25ex]
h_M Bu(t_M)
\end{smallmatrix}\right).
\end{align}
We emphasize that due to dissipativity of $A=(J-R)Q$, the diagonal blocks of the matrix on the left hand side are invertible for any choice of step sizes.
\medskip

\noindent \textbf{Discretized adjoint equation.} 
To discretize the backwards-in-time equation corresponding to $\mathcal{A}^-$, we again project onto the discrete space $X_h$ defined in \eqref{eq:Xh}. That is, we define
$\mathcal{A}^-_h : X_h \to X_h^*$ by
\begin{align}\label{eq:adjdisc}
\begin{split}
\langle \mathcal{A}_h^-\lambda_{h},\varphi_h\rangle_{X_h^*,X_h} := \sum_{i=1}^M h_i \langle -A^\top\lambda_h(t_i), \varphi_h(t_i) \rangle &+ \sum_{i=0}^{M-1} \langle -[\lambda_h]_i, \varphi_h(t_{i})\rangle \\&+ \langle \lambda_h(T), \varphi_h(T)\rangle
\end{split}
\end{align}
for all $\varphi_h\in X_h$. As for the operator $\mathcal{A}_h$, we may decouple the above equations by choosing a test function which vanishes on all but one interval. Thus, this corresponds to the application of an implicit Euler scheme to the backwards-in-time equation modeled by $\mathcal{A}^-$. Correspondingly, the systems to be solved when solving the adjoint equation are governed by the transposed of the matrix in \eqref{eq:adjdisc}, as will be shown in the upcoming section.

\section{Implementation of goal-oriented adaptivity}\label{sec:comp}

\noindent Having defined a variational framework for the formulation of state and adjoint equation, we provide implementation details of the goal-oriented scheme sketched in Section~\ref{sec:ph_goal}. 
In the following, we consider the state $x\in X:=H^1([0,T],\R^n)$ solving the continuous problem~\eqref{eq:state}, i.e., in view of Proposition~\ref{prop:testfunc}, for all $\varphi \in H^1([0,T],\R^n)$
\begin{align*}
\langle \mathcal{A}x,(\varphi,\varphi(0))\rangle_{L^2([0,T],\R^n)\times \R^n} = \int_0^T \langle  Bu(t),\varphi(t)\rangle \,\mathrm{d}t + \langle x^0,\varphi(0)\rangle.
\end{align*}
Furthermore, by $x_h\in X_h$, we denote the solution to the projected problem, i.e.,
\begin{align*}
\langle \mathcal{A}_h x_h,\varphi_h\rangle_{X_h^*,X_h} = \sum_{i=1}^M \int_{\mathbb{I}_i} \langle Bu(t), \varphi(t_i)\rangle\,\mathrm{d}t + \langle x^0,\varphi(t_0)\rangle \quad \forall \varphi \in X_h.
\end{align*}
As a goal for refinement, we consider the quantity of interest $I_\mathrm{loc}(x) $ as defined 
via \eqref{eq:Iloctotal}. Hence, as motivated in \eqref{eq:adjoint}, we introduce an adjoint state $(\lambda,\lambda_0) \in L^2([0,T],\R^n)\times \R^n$ solving
\begin{align}\nonumber 
&\langle \mathcal{A}^*(\lambda,\lambda_0),\varphi \rangle_{H^1([0,T],\R^n)^*,H^1([0,T],\R^n)} 
\\
& \qquad \qquad= \langle I_\mathrm{loc}'(x_h),\varphi\rangle_{H^1([0,T],\R^n)^*,H^1([0,T],\R^n)} \quad \forall \varphi\in H^1([0,T],\R^n). \label{eq:origadj}
\end{align}
As the right-hand-side of this adjoint equation is a functional on $H^1([0,T],\R^n)$ consisting of $L^2$-terms and point evaluations, i.e.,
\begin{align*}
&\langle I'_\mathrm{loc}(x_h),\varphi\rangle_{X^*,X} \\&= 
\changed{2}\sum_{i=0}^{M-1}  \bigg(\int_{t_i}^{t_{i+1}}-u(t)^\top y(t) + 2x(t)^\top \changed{QRQ}x(t)\,\mathrm{d}t + H(x(t_{i+1})) - H(x(t_{i}))\bigg) \times \\
& \bigg(\!\int_{t_i}^{t_{i+1}} \!\!\!\!\!\!u(t) B^\top Q\varphi(t) \!+\! 2x(t)^\top \changed{QRQ} \varphi(t))\,\mathrm{d}t \!+\! H'(x(t_{i+1}))\varphi(t_{i+1}) \!-\! H'(x(t_i))\varphi(t_i)\bigg),
\end{align*}
we may apply Lemma~\ref{lem:higherreg} which implies that $\lambda$ solving \eqref{eq:origadj} is piecewise weakly differentiable. Hence, we consider the projected version as a backwards-in-time equation
\begin{align}\label{eq:projadj}
\langle \mathcal{A}_h^-\lambda,\varphi_h \rangle_{X_h^*,X_h} = \langle I'(x_h),\varphi_h\rangle_{X_h^*,X_h} \qquad \forall \varphi_h \in X_h,
\end{align}
where $\mathcal{A}_h$ is defined in \eqref{eq:adjdisc}. The computation of the adjoint variable, i.e., the projected version of \eqref{eq:adjoint} amounts to solving
\begin{align}\label{eq:system}
\left(\begin{smallmatrix}
I-h_1A^\top & -I & 0 & 0 & \cdots & 0 \\
0 & I-h_2A^\top & -I & 0 & \cdots & 0 \\
0 & 0 & I-h_3A^\top & -I & \cdots & 0 \\
\vdots & \vdots & \ddots & \ddots & \ddots & \vdots \\
\vdots & \vdots & \vdots & 0 & I-h_M A^\top & -I \\
0 & 0 & 0 & 0 & 0 & I \\
\end{smallmatrix}\right)
\begin{pmatrix}
\lambda_0\\ \lambda_1\\ \vdots\\ \lambda_M
\end{pmatrix}
=
\left(\begin{smallmatrix}
\langle I'_\mathrm{loc}(x),\varphi_1\rangle_{X_h^*,X_h}\\
\langle I'_\mathrm{loc}(x),\varphi_2\rangle_{X_h^*,X_h}\\[.2em]
\vdots \\[.7em]
\langle I'_\mathrm{loc}(x),\varphi_M\rangle_{X_h^*,X_h}
\end{smallmatrix}\right).
\end{align}
Note, that the matrix on the left-hand sides is the transposed of the matrix in \eqref{eq:state_disc}, which resembles the adjoint structure of $\mathcal{A}_h$ and $\mathcal{A}_h^-$.

Furthermore, setting $y_i = B^\top Q x_i$ and approximating the integrals by a rectangular rule, we have for all 
$\varphi_h\in X_h$,
\begin{align*}
\langle I'_\mathrm{loc}(x),\varphi_h\rangle_{X_h^*,X_h} &\approx 
\changed{2}\sum_{i=1}^{M-1}  \bigg(h_i(-u_i^\top y_i + x_i^\top \changed{QRQ}x_i) + H(x_{i+1}) - H(x_{i})\bigg) \times
\\
&   (h_i(u_i B^\top Q\varphi_i + 2x_i^\top \changed{QRQ} \varphi_{i}) + H'(x_{i+1})\varphi_{i+1} - H'(x_{i})\varphi_{i}).
\end{align*}

\noindent \textbf{Approximating $\lambda$ from $\lambda_h$.} To evaluate the error estimator as given in \eqref{eq:errorestimator} we have to approximate $\lambda$ solving \eqref{eq:origadj} by means of $\lambda_h$ solving $\eqref{eq:projadj}$. To this end, an efficient strategy, cf.~\cite{MeidVex07}, is to utilize $\lambda \approx I_h^1\lambda_h$, where $I_h^1$ is the piecewise linear interpolation operator. Then,
\begin{align*}
(I_h^1\changed{\lambda_h} - \changed{\lambda_h})(t_i) = \changed{\lambda_h(t_{i+1}) - \lambda_h(t_{i})}.
\end{align*}
Other alternatives include a solution on a refined mesh or a higher-order method, which yield similar results, but are more costly than the linear interpolation technique, see \cite{MeidVex07}.

\noindent \textbf{Evaluation of error estimator.} Approximating the integrals weighted with the interpolant $I_h^1\lambda_h$ by the trapezoidal rule, and the integrals weighted by the piecewise constant function $\lambda_h$ by the right endpoint rule, we obtain the approximation of \eqref{eq:errorestimator} using 
\begin{align*}
& \langle \mathcal{A}(x-x_h),\lambda-\lambda_h\rangle \approx \langle (Bu,x^0) - \mathcal{A}x_h,(I_h^1 - \mathrm{Id})\lambda_h\rangle \\
& =   \sum_{i=1}^M \frac{h_i}{2} \langle -Ax_h(t_i), \lambda_h(t_i) - \lambda_h(t_{i-1})
\rangle +\langle x_h (t_{i}) - x_{h} (t_{i-1}), \lambda_h(t_{i}) - \lambda_h(t_{i-1})\rangle\\
& \qquad + \frac{h_i}{2} \left(\lambda(t_{i-1})^\top Bu(t_{i-1}) - \lambda(t_{i})^\top Bu(t_i) \right).
\end{align*}
Thus, we may define local error indicators via
\begin{align}
\eta_i &:= \Big|\frac{h_i}{2} \langle -Ax_h(t_i), \lambda_h(t_i) - \lambda_h(t_{i-1})
\rangle +\langle x_h (t_{i}) - x_{h} (t_{i-1}), \lambda_h (t_i) - \lambda_h (t_{i-1})\rangle \nonumber \\
& \qquad + \frac{h_i}{2} \left(\lambda(t_{i-1})^\top Bu(t_{i-1}) - \lambda(t_{i})^\top Bu(t_i) \right)\Big|.\label{eq:errorindicators}
\end{align}

\noindent \textbf{Refinement procedure.} To refine the mesh based on the local error indicators above, we use a simple Dörfler criterion \cite[Section 4.2]{Doerfler1996}, which aims to reduce the error by a certain percentage. To this end, we sort the error indicators by their size and refine the corresponding cells until a certain percentage of the error is covered.

We summarize the proposed approach in Algorithm~\ref{alg:goee}.
\begin{algorithm}
	\caption{Goal-oriented adaptivity for pHs}\label{alg:goee}
	\begin{algorithmic}
		\Require Initial time grid $\mathcal{T}$
		\While{refinement needed}
		\State Solve the forward equation \eqref{eq:state_disc} to obtain $x_h$
		\State Solve the backward equation \eqref{eq:system} to obtain $\lambda_h$
		\State Evaluate error indicators $\eta_i$ via \eqref{eq:errorindicators}
		\State Refine intervals corresponding to fixed percentage of error
		\EndWhile
	\end{algorithmic}
\end{algorithm}
There, we observe that the main additional computational effort compared to usual time stepping is the computation of the adjoint state. Thus, in the following section, we provide a parallelizable approximation for the backward equation \eqref{eq:system}.

\section{A dissipativity exploiting parallelizable strategy for the adjoints}\label{sec:para}
\noindent In this section, we provide an efficient way of approximating the sensitivities solving \eqref{eq:system} using the port-Hamiltonian structure.  More precisely, we suggest
\begin{align}\label{eq:system_approx}
\left(\begin{smallmatrix}
I-h_1A^\top &0 & 0 & 0 & \cdots & 0 \\
0 & I-h_2A^\top & 0 & 0 & \cdots & 0 \\
0 & 0 & I-h_3A^\top & 0 & \cdots & 0 \\
\vdots & \vdots & \ddots & \ddots & \ddots & \vdots \\
\vdots & \vdots & \vdots & 0 & I-h_M A^\top & 0 \\
0 & 0 & 0 & 0 & 0 & I \\
\end{smallmatrix}\right)
\begin{pmatrix}
\tilde \lambda_0\\ \tilde \lambda_1\\ \vdots\\ \tilde \lambda_M
\end{pmatrix}
=\left(
\begin{smallmatrix}
\langle I'_\mathrm{loc}(x),\varphi_1\rangle_{X_h^*,X_h}\\
\langle I'_\mathrm{loc}(x),\varphi_2\rangle_{X_h^*,X_h}\\[.2em]
\vdots \\[.7em]
\langle I'_\mathrm{loc}(x),\varphi_M\rangle_{X_h^*,X_h}
\end{smallmatrix}\right)
\end{align}
as an approximation of \eqref{eq:system}.
This can be viewed as the application of one step of a block-Jacobi method initialized at zero and intuitively corresponds to ignoring the error in the quantity of interest on the right-hand side from the next subinterval and only computing local sensitivities. We note that assuming we have factorizations of the matrices $I-h_1A, \ldots, I-h_MA$ at hand (as e.g., computed when solving the forward system \eqref{eq:state_disc}), then solving \eqref{eq:system} amounts to back substitution and hence may be achieved with a complexity of $\mathcal{O}((M+1)n)$. The approximation presented above in \eqref{eq:system_approx}, however, may be completely performed in parallel for all blocks and thus, of course depending on the number of available cores, may be evaluated in $\mathcal{O}(n)$ time. 

Moreover, if we were to apply more iterations of a Jacobi (or Gauss-Seidel method), we get the iteration matrix
\begin{align*}
&G := \left(\begin{smallmatrix}
I-h_1A^\top & 0 & 0 & 0 & \cdots & 0 \\
0 & I-h_2A^\top & 0 & 0 & \cdots & 0 \\
0 & 0 & I-h_3A^\top & 0 & \cdots & 0 \\
\vdots & \vdots & \ddots & \ddots & \ddots & \vdots \\
\vdots & \vdots & \vdots & 0 & I-h_M A^\top & 0 \\
0 & 0 & 0 & 0 & 0 & I \\
\end{smallmatrix}\right)^{-1}
\left(\begin{smallmatrix}
0 & I & 0 & 0 & \cdots & 0 \\[0.5ex]
0 & 0 & I & 0 & \cdots & 0 \\[0.5ex]
0 & 0 & 0 & I & \cdots & 0 \\
\vdots & \vdots & \ddots & \ddots & \ddots & \vdots \\
\vdots & \vdots & \vdots & 0 & 0 & I \\[0.5ex]
0 & 0 & 0 & 0 & 0 & 0 \\
\end{smallmatrix}\right)
\\[0.5ex]
& \hspace{2cm}= 
\left(\begin{smallmatrix}
0 & (I-h_1A^\top)^{-1} & 0 & 0 & \cdots & 0 \\
0 & 0 & (I-h_2A^\top)^{-1} & 0 & \cdots & 0 \\
0 & 0 & 0 & (I-h_3A^\top)^{-1} & \cdots & 0 \\
\vdots & \vdots & \ddots & \ddots & \ddots & \vdots \\
\vdots & \vdots & \vdots & 0 & 0 & (I-h_MA^\top)^{-1} \\
0 & 0 & 0 & 0 & 0 & 0 \\
\end{smallmatrix}\right).
\end{align*}
As this matrix is nilpotent, the corresponding block-Jacobi method for \eqref{eq:system} converges in at most $M+1$ steps. If the matrix-vector products are performed in parallel, this recovers then the complexity $\mathcal{O}((M+1)n)$ of solving the original system \eqref{eq:system} with back substitution.

However, we note that the convergence rate of this method is greatly affected by the stability properties of the underlying port-Hamiltonian system. 
In particular, if $A$ has only eigenvalues with negative real part and with very large absolute value, 
i.e., the system is strongly exponentially stable, the Jacobi method converges with a fast rate. This can be observed by the proof of the following result. 

\begin{proposition}
	If the port-Hamiltonian system is exponentially stable, i.e., if $A = (J-R)Q$ is a Hurwitz matrix, then there exists a matrix norm in which $G$ is a contraction.
\end{proposition}
\begin{proof}
	If $A$ is Hurwitz, its eigenvalues lie in the open left half plane. Correspondingly, for all $h>0$ the eigenvalues of $I-h A^\top = I - h Q(J-R)^\top$ have real part larger than one such that in particular $\min |\mu(I-hA)|>1$, where where $\mu(M)$ denotes the eigenvalues of a matrix $M$. Hence
	\begin{align*}
	\rho((I-hA^\top)^{-1}) = \max |\mu((I-hA^\top)^{-1})| &= \max \frac{1}{|\mu(I-hA)|} = \frac{1}{\min |\mu(I-hA)|}<1.
	\end{align*}
	As for any $\varepsilon> 0$ there is a matrix norm $\|\cdot\|_\varepsilon$ such that 
	\begin{align*}
	\|(I-hA^\top)^{-1}\|_\varepsilon <  \rho((I-hA^\top)^{-1}) + \varepsilon,
	\end{align*}
	this implies that the blocks of $G$ are contractions. Hence, also $G$ is contraction.
\end{proof}
This observation on the discrete level, which estimates the quality of the approximation scheme \eqref{eq:system_approx} in case of exponentially stable dynamics can also be interpreted on the continuous level. To this end, we observe that \eqref{eq:system_approx} includes the same dynamics as \eqref{eq:system}, where however the terminal value is set to zero which corresponds to leaving out the block-super-diagonal of the matrix in \eqref{eq:system} consisting of negative identities.
\begin{proposition}
	Let $\lambda,\tilde \lambda\in H^1([t_i,t_{i+1}],\R^n)$ solve
	\begin{align*}
	- \lambda' &= A^\top \lambda + f \qquad \lambda(t_{i+1}) = \lambda_T\\
	-\tilde \lambda' &= A^\top \tilde\lambda + f \qquad \tilde \lambda(t_{i+1}) = 0.
	\end{align*}
	Then, for all $t\in [t_i,t_{i+1}]$, \changed{and denoting $\|v\|_Q := \sqrt{v^\top Q v}$,}
	\begin{align}\label{eq:normbound}
	\|\lambda(t)-\tilde \lambda(t)\|_Q \leq \|\lambda_T\|_Q, 
	\end{align}
	i.e., the influence of the perturbation is bounded uniformly in the size of the subinterval $h_i = t_{i+1}-t_i$.    If additionally $A$ (or equivalently $A^\top$) is Hurwitz, then there is $\omega> 0$ and $C\geq 1$ such that for all $t\in [t_i,t_{i+1}]$
	\begin{align}\label{eq:errordecay}
	\|\lambda(t)-\tilde \lambda(t)\| \leq Ce^{-\omega (t_{i+1}-t)}\|\lambda_T\|.
	\end{align}
\end{proposition}
\begin{proof}
	The claim directly follows from subtracting the systems for $\lambda$ and $\tilde \lambda$ and applying the variation of constants formula, i.e., for all $t\in (t_{i},t_{i+1}]$
	\begin{align*}
	\lambda(t)-\tilde{\lambda}(t) = e^{(t_{i+1}-t)A}\lambda_T.
	\end{align*}
	As $A = (J-R)Q$ is dissipative in the $\langle \cdot,Q\cdot\rangle$-scalar product (due to skew-symmetry of $J$ and symmetry and positive semi-definiteness of $R$), it generates a contraction semigroup~\cite[Theorem 6.1.7]{JacoZwar12} which implies the claim \eqref{eq:normbound}. Furthermore, if $A$ is exponentially stable, this means that 
	\begin{align*}
	\|\lambda(t)-\tilde{\lambda}(t)\| = \|e^{(t_{i+1}-t)A}\lambda_T\| \leq Ce^{-\omega(t_{i+1}-t)}\|\lambda_T\|
	\end{align*}
	for $\omega = -\max \mathrm{Re}(\mu(A)) > 0$ and $C\geq 1$, where $\mu(A)$ denotes the eigenvalues of $A$.
\end{proof}

\begin{remark}
	Approximating the computation of the sensitivities in a block-diagonal manner is very closely related to DLY (Deuflhard-Leinen-Yserentant; \cite{DeufLein89}) methodologies for elliptic problems. There, due to fast decay of the Green's function in space, the error transport is also localized. Here, the dissipativity of the port-Hamiltonian model in time prevents the global error transport. 
\end{remark}

\section{Numerical results}\label{sec:num}
\noindent In this section, we illustrate our results by means of numerical examples. To this end, we set $n=3$ and $m=1$ and we consider port-Hamiltonian systems given by the matrices
\begin{align}\label{eq:sysmat}
J = \begin{pmatrix}
0&0&1\\
0&0&-1\\
-1&1&0
\end{pmatrix}, \ R_1 = \begin{pmatrix}
1&1&0\\
1&1&0\\
0&0&0
\end{pmatrix}, \ R_2 = I, \ Q = I,\ B = \begin{pmatrix}
1\\1\\0
\end{pmatrix},
\end{align}
i.e., we investigate the case of different dissipation matrices. 
In all experiments, we set the simulation horizon to be $T=10$.

\subsection{Approximation of the sensitivities}
\noindent Here, we illustrate the efficiency of the proposed block-diagonal approximation provided in Section~\ref{sec:para} by means of an example. To this end, we apply the discretization scheme proposed in \eqref{sec:vardisc} to the adjoint equation and solve~\eqref{eq:system} and its approximation \eqref{eq:system_approx} a right-hand side $b$ modeling a local violation of the energy balance at $T/2$, i.e., $b^\top = (0,\ldots,0,5,0,\ldots,0)$. For the system matrix, we compare the results when setting $A = J-R_1$, $A=J-R_2$ and $A = J-10R_2$, cf.~\eqref{eq:sysmat}.

As the eigenvalues of $J-R_1$ are given by $-2, \pm\sqrt{2}i$, the pHs corresponding to this matrix is not exponentially stable due to the presence of two imaginary eigenvalues. Contrary, the eigenvalues of $J-R_2$ are given by $-1,-1\pm\sqrt{2}i$, and the eigenvalues of $J-10R_2$ are given by $-10,-10\pm \sqrt{2}$ such that the underlying port-Hamiltonian systems are exponentially stable.

\begin{figure}[htb]
	\centering
	\includegraphics[width=.32\textwidth]{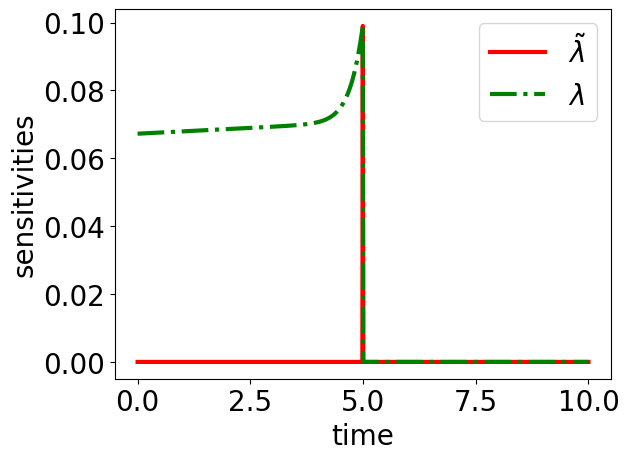}
	\includegraphics[width=.32\textwidth]{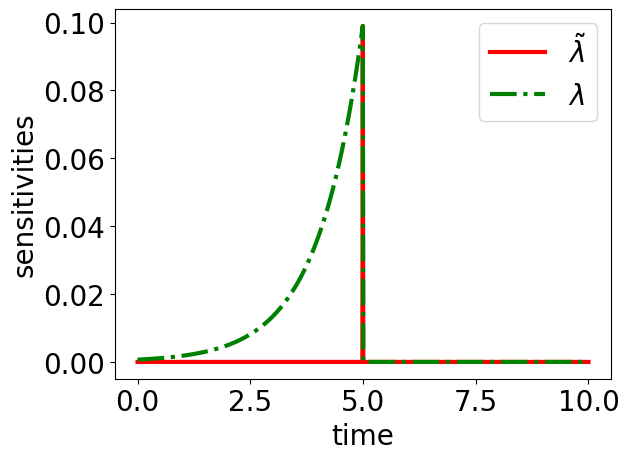}
	\includegraphics[width=.32\textwidth]{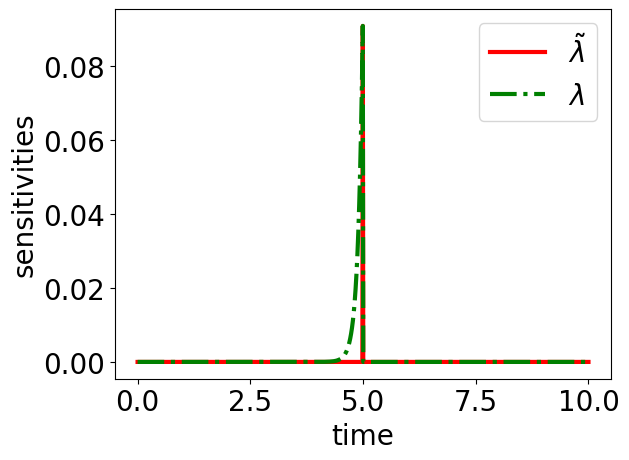}
	\caption{Norm of the sensitivities $\lambda$ and $\tilde{\lambda}$ solving \eqref{eq:system} and \eqref{eq:system_approx} over time. Left: $A=J-R_1$, middle: $A=J-R_2$, right: $A=J-10R_2$.}
	\label{fig:sens}
\end{figure}

In Figure~\ref{fig:sens}, we depict the norm of the corresponding adjoint variables over time for these three different choices of dissipation matrices. We observe in the left plot that for $A=J-R_1$, the perturbation when solving~\eqref{eq:system} is propagated to the initial part of the time horizon. This is due to the two purely complex eigenvalues of $J-R_1$ implying that also errors can be transported without damping. However, when solving the block-diagonal approximation \eqref{eq:system_approx} the influence of this local perturbation is only local due to the decoupled system. 
The  approximation error, i.e., the difference of $\tilde \lambda$ and $\lambda$ however is bounded due to \eqref{eq:normbound}. In the middle plot of Figure~\ref{fig:sens} we see the results for the exponentially stable dynamics governed by $A=J-R_2$. Correspondingly, and in view of \eqref{eq:errordecay}, the error decreases. When increasing this exponential stability by considering $J-10R_2$, the error decay is stronger, as to be expected. 

\subsection{Performance of the goal-oriented scheme: Toy example}
\noindent We compare the results obtained by goal-oriented adaptivity with a uniform refinement strategy and the \textsc{scipy} solver \textsc{odeint}. For the latter, we can obtain a numerical solution on a given time grid, where however, internally more time steps are used. The time stepping is performed using either a BDF-method or an Adams method~\cite{Petz83}. 
To this end, we consider the system given by $A=J-R_1$, that is, only the first dissipation matrix such that the system is not exponentially stable. We choose the initial value $x^0 = \begin{pmatrix}
1&2&1
\end{pmatrix}^\top$ and the input function
\begin{align*}
u(t) = \begin{cases}
0& t\in [0,5]\\
10& t\in (5,10].
\end{cases}
\end{align*}
In the following, we evaluate the performance goal-oriented refinements using Algorithm~\ref{alg:goee} and when choosing the QoI $I_\mathrm{loc}$ defined in \eqref{eq:Iloc}, the weighted QoI $I_\mathrm{loc,\rho}$ defined in \eqref{eq:Ilocrho} and sensitivities computed via \eqref{eq:system}, as well as sensitivities computed via the block-Jacobi-type approximation \eqref{eq:system_approx}.

\begin{figure}[H]
	\centering
	\includegraphics[width=.48\linewidth]{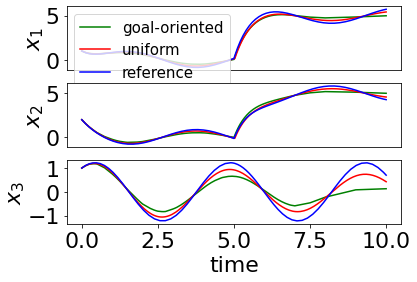}\hspace{.03\linewidth}
	\includegraphics[width=.48\linewidth]{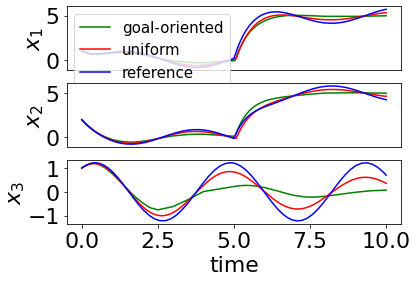}\\
	\includegraphics[width=.48\linewidth]{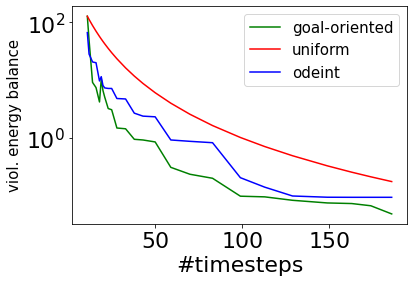}\hspace{.03\linewidth}
	\includegraphics[width=.48\linewidth]{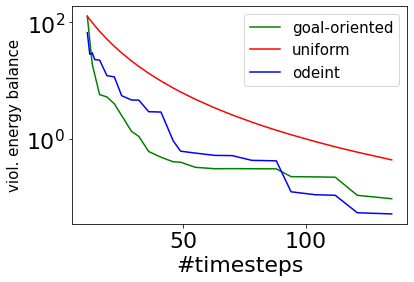}
	\caption{Trajectories (top) and violation of energy balance (bottom) with QoI $I_{\mathrm{loc}}$. Left column: sensitivities computed form \eqref{eq:system}, right column: sensitivities computed from \eqref{eq:system_approx}}
	\label{fig:1}
\end{figure}

First, in Figure~\ref{fig:1}, we evaluate the performance of $I_\mathrm{loc}$ with sensitivities computed from \eqref{eq:system} with a scheme using approximated sensitivities from \eqref{eq:system_approx}.

\noindent For the sensitivities computed via the exact adjoint system \eqref{eq:system}, we observe in the bottom left plot that the violation of the energy balance is always lowest when using the goal-oriented refinement regime. We stress that \textsc{odeint} internally utilizes a higher number of grid points. In the upper plots of Figure~\ref{fig:1} we observe that despite this good performance in view of the goal of refinement, the error of the scheme is relatively large in the third component towards the end of the horizon, as the refinement w.r.t.~$I_\mathrm{loc}$ solely focuses on the error in the energy balance. In the right column of Figure~\ref{fig:1}, we plot the results using approximated sensitivities. Whereas for a smaller amount of time steps, the energy balance violation is still the lowest, we reach a point of saturation compared to \textsc{odeint}. Intuitively, this indicates that the dissipation is not strong enough to justify localizations corresponding to the block-diagonal approximation. We note that the considered matrix $A=(J-R_1)Q$ is not Hurwitz due to the presence of two imaginary eigenvalues, that is, the system is not exponentially stable. 

\begin{figure}[htb]
	\centering
	\includegraphics[width=.48\linewidth]{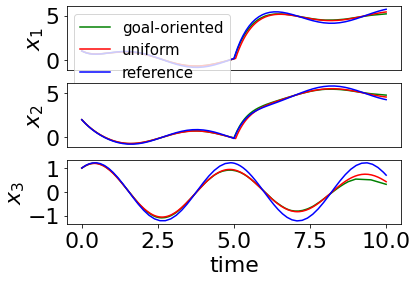}\hspace{.03\linewidth}
	\includegraphics[width=.48\linewidth]{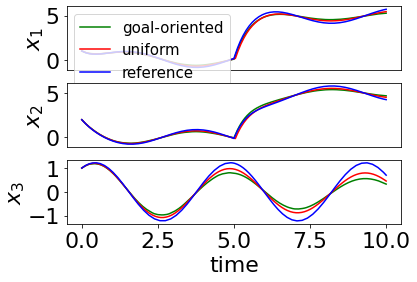}\\
	\includegraphics[width=.48\linewidth]{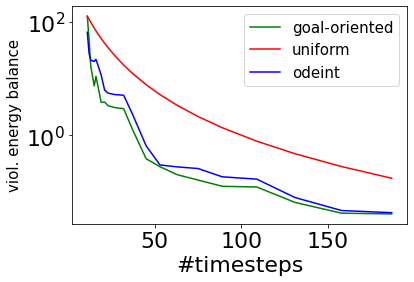}\hspace{.03\linewidth}
	\includegraphics[width=.48\linewidth]{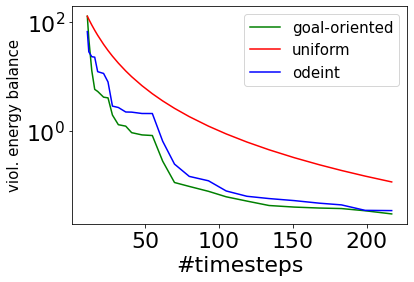}
	\caption{Trajectories (top) and violation of energy balance (top) with QoI $I_{\mathrm{loc},\rho}$. Left column: sensitivities computed form \eqref{eq:system}, right column: sensitivities computed from \eqref{eq:system_approx}}
	\label{fig:2}
\end{figure}

\noindent In Figure~\ref{fig:2}, we depict the same quantities as in Figure~\ref{fig:1}, where however we use the weighted quantity of interest $I_\mathrm{loc,\rho}$ with $\rho=10$. As to be expected, in e.g., in the left column, we observe a better performance of the goal-oriented method in terms of the trajectory error compared to Figure~\ref{fig:1}. However, this better trajectory error comes at the cost of a higher energy balance violation.

\begin{figure}[htb]
	\centering
	\includegraphics[width=.48\linewidth]{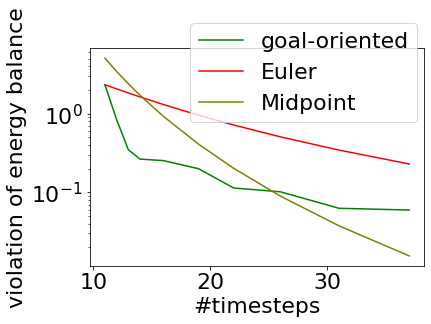}\hspace{.03\linewidth}
	\includegraphics[width=.48\linewidth]{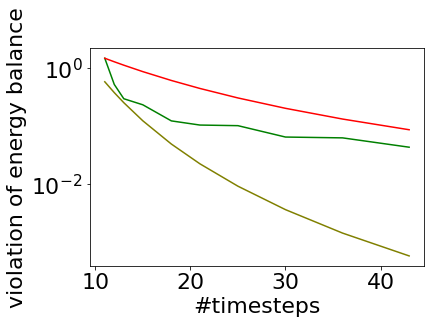}
	\caption{Violation of energy balance and comparison to implicit midpoint rule with uniform refinement with no control input: Left: $R$ from \eqref{eq:sysmat}, Right: Dissipation given by $0.6\cdot R$.}
	\label{fig:comp}
\end{figure}

\noindent \changed{Last, we also provide a comparison to the implicit midpoint rule which is known to be one of the simplest symplectic and (due to the quadratic Hamiltonian) energy-preserving integrators~\cite{HairHoch06,MehrMora19,KotyLefe19}. We consider the case without input, that is $u\equiv 0$ and compare the performance for different dissipation matrices. In case of no dissipation, the implicit midpoint rule leads to no violation of the energy balance. In Figure~\ref{fig:comp} we depict the results. On the left we see that, for the dissipation matrix considered earlier, the implicit midpoint rule overtakes the goal-oriented adaptive implicit Euler scheme. This is due to the implicit midpoint rule enjoying a higher convergence rate by one order. On the right-hand side however, we see that for smaller dissipation, the implicit midpoint rule outperforms the goal-oriented scheme for all choices of time steps. This clearly motivates future research formulating goal-oriented adaptive midpoint schemes: The adaptivity then particularly targets time instances at which energy is dissipated (or supplied) and thus complements the symplectic properties of the standard implicit midpoint scheme.}

\allowdisplaybreaks
\subsection{Performance of the goal-oriented scheme: Electrical circuit}
\noindent \changed{In this part, we consider a circuit model in the sense of a modified nodal analysis as depicted in Figure~\ref{fig:circuit}.}
\begin{figure}[htb]
	\hspace{-1.2cm}
	\begin{tikzpicture}
	\begin{circuitikz}[scale=.7,transform shape]
	%
	%
	\draw (-0.25,3) to[I, l_=$\imath(t)$] (-0.25,0);
	\draw (-0.25,3) -- (2.,3);    
	\draw[fill=black] (2,3) ellipse (.05 and .05) node[above,yshift=0.1cm]{$e_{0}$};
	\draw (1.6,2.5) to[C, l_=$C$] (1.6,0.5);
	\draw (2.4,2.5) to[R, l=$R_1$] (2.4,.5);
	\draw (1.6,2.5) -- (2.4,2.5);
	\draw (1.6,0.5) -- (2.4,0.5);
	\draw (2,0.5) -- (2, 0);
	\draw (2,2.5) -- (2, 3);
	
	\draw (2.,3) to[R, l=$R_d$] (4.5,3);
	\draw[fill=black] (4.5,3) ellipse (.05 and .05) node[above,yshift=0.1cm]{$e_{1}$};
	\draw (4.5,3) to[C, l=$C_R$] (4.5,0);
	\draw (4.5,3) to[L, l=$L$] (7.,3);
	\draw[fill=black] (7,3) ellipse (.05 and .05) node[above,yshift=0.1cm]{$e_{2}$};
	
	\begin{scope}[xshift=5cm]
	\draw (1.6,2.5) to[C, l_=$C$] (1.6,0.5);
	\draw (2.4,2.5) to[R, l=$R_1$] (2.4,.5);
	\draw (1.6,2.5) -- (2.4,2.5);
	\draw (1.6,0.5) -- (2.4,0.5);
	\draw (2,0.5) -- (2, 0);
	\draw (2,2.5) -- (2, 3);
	\end{scope}
	
	\draw (7.,3) -- (8.,3);
	\draw (-0.25,0) -- (8.,0);
	\draw[fill=black] (8.5,3) ellipse (.05 and .05);
	\draw[fill=black] (8.75,3) ellipse (.05 and .05);
	\draw[fill=black] (9.,3) ellipse (.05 and .05);
	\draw[fill=black] (8.5,0) ellipse (.05 and .05);
	\draw[fill=black] (8.75,0) ellipse (.05 and .05);
	\draw[fill=black] (9.,0) ellipse (.05 and .05);

	\begin{scope}[xshift=7.5cm]
	\draw (2.5,3) to[R, l=$R_d$] (5,3);
	\draw[fill=black] (5,3) ellipse (.05 and .05) node[above,yshift=0.1cm]{$e_{2n-1}$};
	\draw (5,3) to[C, l=$\!C_R$] (5,0);
	\draw (5,3) to[L, l=$L$] (7.5,3);
	\draw[fill=black] (7.5,3) ellipse (.05 and .05) node[above,yshift=0.1cm]{$e_{2n}$};
	
	\draw (7.5,3) -- (9.5,3);
	\draw (2.5,0) -- (9.5,0);
	
	\begin{scope}[xshift=5.5cm]
	\draw (1.6,2.5) to[C, l_=$C$] (1.6,0.5);
	\draw (2.4,2.5) to[R, l=$R_1$] (2.4,.5);
	\draw (1.6,2.5) -- (2.4,2.5);
	\draw (1.6,0.5) -- (2.4,0.5);
	\draw (2,0.5) -- (2, 0);
	\draw (2,2.5) -- (2, 3);
	\end{scope}
	

	\draw (9.5,3) to[R, l=$R_L$] (9.5,0);
	
	\end{scope}
	
	\begin{scope}[xshift=-1cm]
	\draw (3.,0) -- (3.,-0.6);
	\draw (2.6,-.6) -- (3.4,-.6);
	\draw (2.75,-.7) -- (3.25,-.7);
	\draw (2.9,-.8) -- (3.1,-.8);
	\draw[white] (9,5) ellipse (.01 and .01);
	\end{scope}
	
	\draw[decorate,decoration={brace, mirror}] (2.6,-0.3) -- (7.5,-.3) node[midway,below,yshift=-0.1cm]{$n$ times};
	%
	\end{circuitikz}
	\end{tikzpicture}
	\label{fig:circuit}
	\caption{Depiction of the considered electrical circuit}
\end{figure}
\changed{The variables are the node potential $e_i$ ($i=0,1,\dotsc, 2n$), and the currents through the inductors $\jmath_i$ ($i=1,\dotsc, n$). 
	The corresponding pHs formulation is given by 
	\begin{align*}
	E\dot{x}(t) = (J-R) x (t)+ Bu(t)
	\end{align*}
	with the state variable
	\[x = [e_0, e_1, \jmath_1, e_2,e_3,\jmath_2, \dotsc, e_{2n-1},\jmath_n, e_{2n}]^\top
	\in \real^{3n+1}
	\]
	the (invertible) descriptor matrix
	\[
	E=\text{blkdiag}(C,E_0, \dotsc,E_0) \in \real^{(3n+1)\times (3n+1)}\ \ \text{with} \ \
	E_0 = \text{diag} (C_R, L, C) \in \real^{3\times 3},
	\]
	the structure matrix
	\[
	J=\text{blkdiag}(0,J_0, \dotsc,J_0) \in \real^{(3n+1)\times (3n+1)}\ \ \text{with} \ \
	J_0 = \left(\begin{smallmatrix} 0 & -1 & 0 \\
	1 &  0 & -1 \\
	0 &  1 &  0
	\end{smallmatrix}\right),
	\]
	the dissipation matrix
	\[
	R = \text{diag}(R_0, \dotsc, R_0, \tfrac{1}{R_L}) \in \real^{(3n+1)\times (3n+1)} \ \ \text{with} \ \
	R_0 = \left(\begin{smallmatrix}  \tfrac{1}{R_1}+\tfrac{1}{R_d} & -\tfrac{1}{R_d} & 0 \\
	-\tfrac{1}{R_d} &  \tfrac{1}{R_d} & 0 \\
	0 &  0 & 0
	\end{smallmatrix}\right)
	\]
	and the input matrix corresponding to a voltage input $u = \imath$
	\[
	B= \begin{pmatrix} 1&0&\dotsc& 0\end{pmatrix}^\top.
	\]
	We inspect the time interval $[0,\,  5\cdot 10^{-6}]$, and consider the parameters $C=5 \cdot 10^{-7} , \; C_R=5\cdot 10^{-9}, \; L=10^{-6},\;  R_d=5,\; R_1=200,\; R_L=10,$ and
	\[
	\imath(t) = \sin(5\cdot 10^{6} \pi t), \; x(0) = 10^{-7}\changedsecond{\begin{pmatrix} 1 & \ldots & 1 \end{pmatrix}^\top}.
	\]
	Furthermore, we choose $n=5$ leading to a 16-dimensional system, however we stress that the results are very robust with respect to adding more elements to the circuit.}

\begin{figure}[htb]
	\centering
	\includegraphics[width=0.48\linewidth]{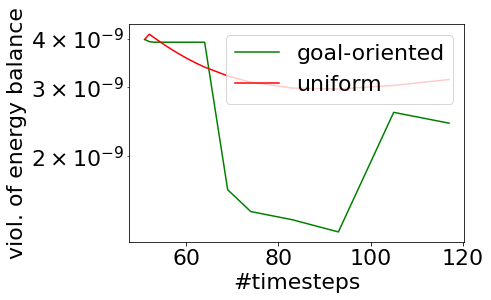}\hspace{0.03cm}
	\includegraphics[width=0.48\linewidth]{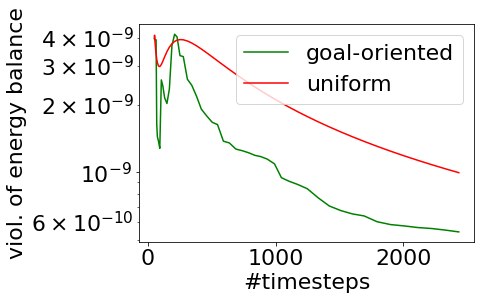}\\
	\includegraphics[width=0.48\linewidth]{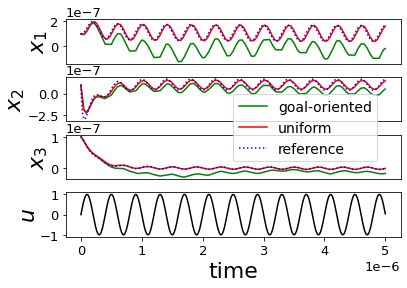}\hspace{0.03cm}
	\includegraphics[width=0.48\linewidth]{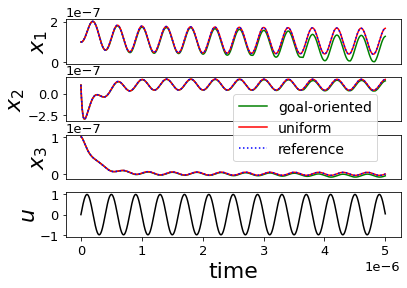}
	\caption{Computational results for the electrical circuit. Left: 20 refinement loops; Right: 60 refinement loops. The vertical axis of the top figures is normalized in view of the total energy supplied to the system.}
	\label{fig:circ_res}
\end{figure}

\changed{We depict the results in Figure~\ref{fig:circ_res}. In the left column, we utilize 20 refinement cycles and observe that the goal-oriented scheme outperforms the uniform refinements in view of the energy balance. However, in view of the trajectories over time, there is clearly a mismatch in the first coordinate. This mismatch is \changedsecond{greatly reduced} when performing more refinement cycles as depicted on the right. 
	The reference trajectory is computed with a high-order BDF scheme on a very fine grid using 10001 time steps.}

\section{Conclusion and outlook}
\noindent We proposed an adaptive and goal-oriented grid refinement approach tailored to the energy balance in port-Hamiltonian systems. To this end, we formulated the state equation and the corresponding adjoint equation as variational problems and proved a result for higher regularity of the adjoint for right-hand sides occurring in the proposed scheme. To guarantee an efficient evaluation of the error estimator, we further provided a block-diagonal and stability-exploiting approximation for the adjoint system which we motivated from a linear algebraic perspective via a block-Jacobi method. Then, we provided a numerical example highlighting that the proposed scheme leads to smaller violations of the energy balanced when compared to step size controlled or higher-order time stepping schemes. 

In the future, we plan to consider goal-oriented multirate schemes~\cite{BaSc_Estep2012} for pHs and to relate the dissipativity-based arguments of Section~\ref{sec:comp} with the approach using generalized Green's functions to capture error transport as proposed in \cite{EsteHols05}. \changed{Furthermore, a combination with a symplectic integrator is expected to lead to a symbiosis of grid adaptivity and geometric integration.}
\bibliographystyle{abbrv}
\bibliography{references}
\end{document}